\documentclass[11pt]{article}

 \usepackage{amsfonts, amsthm, amssymb }
\usepackage[dvips]{epsfig,graphics}

 \newtheorem{theorem}{Theorem}[section]
 \newtheorem{corollary}[theorem]{Corollary}
 
 \newtheorem{proposition}[theorem]{Proposition}
\theoremstyle{definition}

 \newtheorem{remark}[theorem]{Remark}

\def\l{\langle}
\def\r{\rangle}
\def\S{\mathbb S}
\def\Z{\mathbb Z}
\def\D{\mathbb D}
\newcommand{\Fig}{Fig}

\begin{document}

\title{When does a satellite knot fiber? }
\author{
Mikami Hirasawa\footnote{partially supported by MEXT, Grant-in-Aid for
Young Scientists (B) 18740035}\\
Department of Mathematics\\ Nagoya Institute of Technology\\Nagoya Aichi 466-8555\\ JAPAN\\
\and
Kunio Murasugi\footnote{partially supported by NSERC Grant A 4034}\\
Department of Mathematics\\ University of Toronto\\ Toronto ON M5S~2E4 \\CANADA
\and
Daniel S. Silver\footnote{Supported in part by NSF Grant DMS \#0304971}\\ Department of Mathematics and Statistics
\\University of South Alabama\\ Mobile AL  36688\\ USA
}

\maketitle

\begin{abstract}

Necessary and sufficient conditions are given 
for a satellite knot to be fibered. 
Any knot $\tilde k$ embeds in an unknotted solid torus 
$\tilde V$ with arbitrary winding number 
in such a way that no satellite knot with pattern 
$(\tilde V, \tilde k)$ is fibered. 
In particular, there exist nonfibered satellite knots 
with fibered pattern and companion knots and nonzero winding number. 
\end{abstract}



\section{Introduction} Proposition 4.15 of the monograph {\it Knots} 
\cite{bz} erroneously asserts and proves that a satellite knot is fibered if and only if it has nonzero winding number and both  its companion and pattern knots are fibered. In 1990, the third author of the present paper contacted  the authors of {\it Knots}, informed them of the mistake,  and provided a counterexample. The recent second edition of the now-classic textbook acknowledged the error and reproduced the example. However, the  necessary and sufficient conditions that it provided are not the most convenient in practice. 

The present paper is an expanded version of an unpublished note that the third author wrote in 1991. Its algebraic necessary and sufficient conditions for a satellite knot to fiber, do not appear to be well known. (It is briefly mentioned on pages 221--222 of \cite{silver}.)

Theorem 4.2 of the monograph \cite{en} by D. Eisenbud and W. Neumann gives topological necessary and sufficient conditions for a satellite knot in a homology $3$-sphere to fiber. Stated in terms of ``splice components" arising from the Jaco-Shalen Johannson splitting theorem, it requires that each component be fibered and that the fibers meet each other compatibly. 

An alternative, topological condition presented here (see Theorem \ref{alternate}) is more convenient in many situations. It was discovered by the first author, who proved it using D. Gabai's sutured manifold theory. An elementary group theoretic proof was subsequently found by the second author.

We work throughout in the smooth category. As usual, a {\it knot} is an embedded circle in the $3$-sphere ${\mathbb S}^3$.  It is {\it trivial} if it is the boundary of an embedded disk in ${\mathbb S}^3$.
More generally, a {\it link} is a union of  finitely many pairwise disjoint knots. When each component is provided with a direction, the link is said to be {\it  oriented}.

The complement of an open  tubular neighborhood of a knot or  link $\ell$ is called its {\it exterior}, and it is denoted by $X(\ell)$. Its boundary $\partial X(\ell)$ can be identified with $\ell \times {\Bbb S}^1$. A knot $k$ is {\it fibered} if the projection $\partial X(k) \to {\Bbb S}^1$ extends to a fibration of $X(k)$. In such a case, $X(k)$ is the total space of a fiber bundle over the circle $\S^1$. By theorems of L. Neuwirth \cite{neuwirth} and J. Stallings \cite{stallings}, $k$ is fibered if and only if the commutator subgroup $[\pi_k, \pi_k]$ of the knot group $\pi_k= \pi_1(\S^3\setminus k)$ is finitely generated and hence free. 
The  commutator subgroup is the kernel of the abelianization homomorphism 
$\chi: \pi_k \to \Z$ sending the class of an oriented meridian of $k$ to $1$.


\section{Fibering criterion.} Satellite knots were introduced by H. Schubert \cite{schubert} as a generalization of knot product. We briefly review their definition and some key properties. Let $\tilde k$ be a knot embedded in a standard solid torus $\tilde V  = {\S}^1 \times {\mathbb D}^2 \subset {\S}^3$. Assume that $\tilde k$ is not isotopic to $\S^1 \times \{0\}$ nor is it contained in any $3$-ball of $\tilde V$. Let $h$ be a homeomorphism from $\tilde V$ onto a closed tubular neighborhood of a nontrivial knot $\hat k$, mapping a longitude of $\tilde V$ onto a longitude of $\hat k$. (By a {\it longitude} of $\tilde V$, we mean a meridian of the complementary solid torus ${\S}^3 \setminus int\ \tilde V$.)
The image $k=h(\tilde k)$ is called a {\it satellite knot} with {\it companion knot} $\hat k$ and {\it pattern} $(\tilde V, \tilde k)$. By a
{\it pattern knot} we mean $\tilde k$ regarded in ${\S}^3$. 
The {\it winding number} of $\tilde k$ in $\tilde V$ is the nonnegative integer $n$ such that the homomorphism $H_1 \tilde k  \to H_1 \tilde V \cong {\Z}$ induced by inclusion has image $n{\Z}$. 

Two patterns $(\tilde V, \tilde k_1)$ and $(\tilde V, \tilde k_2)$ are regarded as the same if $\tilde k_1$ and $\tilde k_2$ are isotopic in $\tilde V$. 

We denote the various groups $\pi_1({\S}^3 \setminus k), 
\pi_1({\S}^3 \setminus \hat k), \pi_1({\S}^3 \setminus \tilde k)$ and $\pi_1(\tilde V\setminus \tilde k)$ by $G, \hat G, \tilde G$ and $H$, respectively.  By the Seifert-van Kampen theorem,

\begin{equation}G\cong H*_{\pi_1(\partial \tilde V)}\hat G,\label{G}\end{equation}
the amalgamated free product of $H$ and $\hat G$ with subgroups
$\pi_1(\partial \tilde V) \subset H$ and $\pi_1(\partial X(\hat k)) \subset \hat G$ identified via $h$. 

It is convenient to regard $H$ as the group of the $2$-component link $\tilde k \cup \tilde m$, where $\tilde m$ is a meridian of $\tilde V$. Let $\tilde l$ be a longitude of $\tilde V$, and let $\lambda$ denote its homotopy class in $H$. The quotient group $\l H \mid \lambda \r$ is naturally isomorphic to $\tilde G$. We define ${\cal K}$ to be the kernel of $\phi\circ \pi$, where $\pi: H \to \l H \mid \lambda \r$ and $\phi: \tilde G \to {\Z}$ are the quotient and abelianization homomorphisms, respectively. When the winding number $n$ is nonzero, then the Reidemeister-Schreier method can be used to show that the commutator subgroup $G'$ decomposes as an amalgamated free product:
\begin{equation}G' \cong {\cal K} *_F (\hat G' * t \hat G' t^{-1} * \cdots * t^{n-1} \hat G' t^{-(n-1)})\label{G'}\end{equation}
Here $t \in G$ is the class of a meridian of $k$, while $F$ is the subgroup of ${\cal K}$ generated by \break $\lambda, t \lambda t^{-1}, \ldots, t^{n-1}\lambda t^{-(n-1)}$. (For a proof, see Section 4.12 of \cite{bz}.)

Let $Y(\tilde k)$ denote the complement of an open tubular neighborhood of $\tilde k$ in $\tilde V$. We say that the pattern $(\tilde V, \tilde k)$ is {\it fibered} if there exists a fibration $Y(\tilde k) \to \S^1$ inducing the homomorphism $\phi\circ \pi$ on fundamental groups. It is well known (and not difficult to see) that in this case, the winding number of $\tilde k$ in $\tilde V$ is nonzero. Moreover, ${\cal K}$ is the fundamental group of the fiber, a compact surface with boundary, and hence ${\cal K}$ is free.

The following fibering criterion for satellite knots is implicit in the knot theory literature (e.g. \cite{en}), but to our knowledge it has not previously appeared. 

\begin{theorem}\label{criterion1} 
Let $k$ be a satellite knot with companion knot $\hat k$ and pattern $(\tilde V, \tilde k)$. 
Then the following conditions are equivalent. 
\begin{enumerate}
\item\quad  $k$ is fibered; 
\item \quad  $\hat G'$ and ${\cal K}$ are finitely generated;
\item \quad  $\hat k$ and $(\tilde V, \tilde k)$ are fibered. 
\end{enumerate}
\end{theorem}

{\it Proof.} Suppose that $k$ is fibered. The winding number $n$ of $\tilde k$ in  $\tilde V$ must be nonzero. For otherwise, $\pi_1(\partial \tilde V) \cong \Z\oplus \Z$, which is a subgroup of $G$ by equation \ref{G}, would be contained in the free commutator subgroup $G'$ -- an impossibility. Consequently, $G'$ has the form 
described by equation \ref{G'}. Since a free product of groups amalgamated over a finitely generated group is itself finitely generated if and only if each of its factors is (see Lemma 4.14 of \cite{bz}, for example), both $\hat G'$ and ${\cal K}$ are finitely generated. 

If $\hat G'$ and ${\cal K}$ are finitely generated, then Stallings's theorem \cite{stallings} implies that $\hat k$ and $(\tilde V, \tilde k)$ are fibered. 

Finally, if $\hat k$ and $(\tilde V, \tilde k)$ are fibered, then $\hat G'$ and ${\cal K}$ are isomorphic to fundamental groups of fibers for respective fibrations. Consequently,  both groups are finitely generated. Equation \ref{G'} implies that  the commutator subgroup $G'$ is finitely generated. By Stallings's theorem \cite{stallings}, the satellite knot $k$ is  fibered. \qed

\section{Pattern fibering criterion.} Let $\ell= \ell_1 \cup \cdots \cup \ell_\mu$ be an oriented link. A meridian of any component acquires an orientation by the right-hand rule. The link is {\it fibered} if its exterior $X(\ell)$ admits a fibration $X(\ell) \to \S^1$. Such a fibration induces an epimorphism $H_1 X(\ell) \to \Z$, and as usual we require that the class of each meridian be sent to $1$. 

The epimorphism $H_1 X(\ell) \to \Z$ is the abelianization of an epimorphism $\chi: \pi_1(\S^3 \setminus \ell) \to \Z$. By \cite{neuwirth} and \cite{stallings}, the link $\ell$ is fibered if and only if the kernel ${\cal A}$ of $\chi$ is finitely generated. The group ${\cal A}$ is known as the {\it augmentation subgroup} of the group of $\ell$. 

\begin{remark}\label{oriented} Orientation is not required in the definition of a fibered knot, since any knot is fibered if and only if the knot with reversed orientation is fibered. Similarly, an orientated link is fibered if and only if the link obtained by reversing the orientation of every component is fibered.
However, reversing the orientation of some but not all components can destroy the condition of fiberedness. (A simple example is provided by the closure of the 
$2$-braid $\sigma_1^4$.) \end{remark}

Consider a pattern $(\tilde V, \tilde  k)$ such that $\tilde k$ has nonzero winding number $n$. We associate an oriented $3$-component link $\tilde k \cup \tilde m \cup \tilde m'$, where $\tilde m$ and $\tilde m'$ are disjoint meridians of $\tilde V$ with opposite orientations. The orientation of $\tilde k$ is arbitrary. 

\begin{theorem}\label{alternate} The pattern $(\tilde V, \tilde k)$ is fibered if and only if $\tilde k \cup\tilde m \cup \tilde m'$ is a fibered link. \end{theorem}

{\it Proof.} As in section 2, let $H$ be the group of the $2$-component link $\tilde k \cup \tilde m$, and let ${\cal K}$ be the kernel of the epimorphism $H \to\Z$ mapping the class of a meridian of $\tilde k$ to $1$ and $\lambda$, the class of a meridian of $\tilde m$, to $0$. 

It suffices to show that the augmentation subgroup ${\cal A}$ of the link $\tilde k \cup\ \tilde  m\  \cup\  \tilde m'$ is the free product of ${\cal K}$ and a finitely generated group (in fact, a free group of rank $n$, the winding number of $\tilde k$ in $\tilde V$). Since the free product of two groups is finitely generated if and only if each of its factors is finitely generated, it will then follow that ${\cal K}$ is finitely generated if and only if ${\cal A}$ is. Equivalently, the pattern $(\tilde V, \tilde k)$ is fibered if and only if the link $\tilde k \cup\tilde m \cup \tilde m'$ is fibered. 

\begin{figure}
\begin{center}
\includegraphics[height=2.3in]{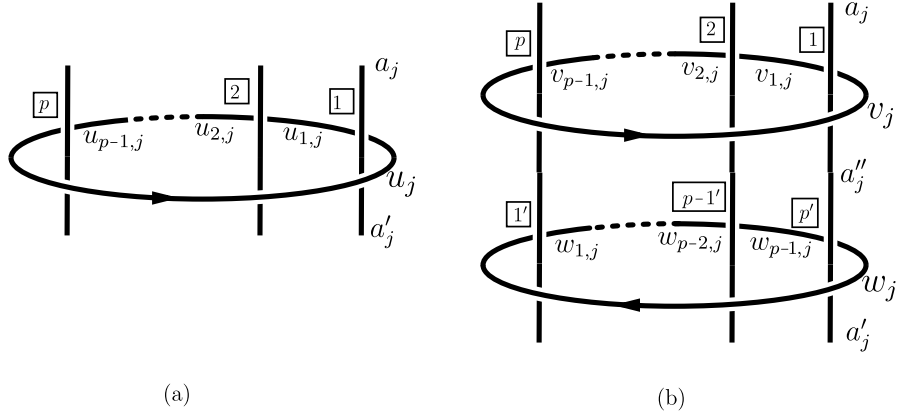}
\caption{(a) Partial diagram of ${\cal D}_1$,
 (b) Partial diagram of ${\cal D}_2$.}
\label{D2}
\end{center}
\end{figure}

Part of a diagram ${\cal D}_1$ for $\tilde k\  \cup  \tilde m$ appears in 
\Fig. \ref{D2} (a) 
with families of generators indicated. There are $p \ge n$ strands of $\tilde k$ passing 
the meridional disk bounded by the unknotted circle representing $\tilde m$.
Assign a weight $\omega$ of $1$ or $0$ to each arc of ${\cal D}_1$ according to whether the arc belongs to $\tilde k$ or $\tilde m$, respectively. 
The Wirtinger algorithm combined with the Reidemeister-Schreier method enables us to write a presentation for ${\cal K}$. Each arc represents an infinite family $a_j, a'_j, b_j, b'_j \ldots, u_j, $  $u_{1, j}, \ldots, u_{p -1, j}$ of generators, indexed by $j \in \Z$. The letters $a, a', b, b',  \cdots$ correspond to arcs of $\tilde k$, with one arc labeled only with the identity element; the remaining letters correspond to arcs of $\tilde m$.

Relators, which come in families indexed by $j\in \Z$, correspond to the crossings of ${\cal D}_1$, and have the form 
$w_j y_{j+\omega} = z_j w_{j+\omega'}$, where $\omega$ and $\omega'$ are the weights of the arcs corresponding to $w_j$ and $z_j$, respectively. We use relators $r_1, \ldots, r_{p -1}$ (indicated in the figure by boxed numbers) to eliminate $u_{1, j}, \ldots, u_{p -1, j}$. 

Recall that in a link diagram, any one Wirtinger relator is a consequence of the others. Consequently, any single relator family is a consequence of the others. We regard $r_p$ as a redundant family of relators.  

The relator families arising from the remaining crossings in  the figure have the form $a'_j = u_j a_j \bar u_{j+1}, b'_j = u_j b_j \bar u_{j+1}, \ldots$, where $\bar{}$ denotes inverse. 

Consider now the partial diagram ${\cal D}_2$ for $\tilde k\ \cup \tilde m\ \cup\ \tilde m'$ that appears in \Fig. \ref{D2} (b), 
again with families of generators indicated. Letters $a, a', a'', b, b', b'', \ldots$ correspond to arcs of $\tilde k$ while $v, v_1, \ldots, v_{p -1}$ and $w, w_1, \ldots, w_{p -1}$ correspond to arcs of $\tilde m$ and $\tilde m'$, respectively.
The parts  of ${\cal D}_1$ and ${\cal D}_2$ not shown are identical.

Relator families $r_1, \ldots, r_{p -1}$ (indicated in the figure by boxed numbers) are used to eliminate $v_{1,j}, \ldots, v_{p -1, j}$. The relator family $r_p$ expresses $v_{j+n}$ in terms of $v_j, \ldots, v_{j+n -1}, a_j, b_j, \ldots$; it also expresses $v _j$ in terms of $v_{j+1}, \ldots, v_{j+n}, a_j, b_j, \ldots$. 

The top front relator families are used to eliminate $a''_j, b''_j, \ldots$. 
Relator families $r_1', \ldots, r_{p-1}'$  are used to eliminate $w_{1, j}, \ldots, w_{p -1, j}$. The  relator family $r_p'$ is  regarded as redundant. 

The bottom front relator families are of the form
$a'_j = \bar w_{j-1} v_{j-1} a_j \bar v_j w_j, b'_j = \bar w_{j-1} v_{j-1} b_j \bar v_j w_j, \ldots$.

Introduce a new generator  family $u_j$ and  defining relators
$u_j = \bar w_{j-1} v_{j-1}$. Rewrite the relators as $w_j = v_j \bar u_{j+1}$, and use them to eliminate $w_j$. 

Relator family $r_p$ can be used to eliminate $v_j$ for  $j  \ne 0, 1,2,\ldots, n-1$, writing them in terms of $v_0, v_1, \ldots, v_{n-1}, a_j, b_j, \ldots$. None of the  remaining relators contain $v_0, \ldots, v_{n -1}$, and so these  generate a free factor $F_n$ of rank $n$.  

The other generators and remaining relators are identical with those of the presentation obtained for ${\cal K}$. Hence ${\cal A}$ is isomorphic to the free product ${\cal K}*F_n$. 
\qed


\section{Examples.} 
We give examples to demonstrate the applicability
of Theorem \ref{alternate}.
We show that the knot $7_6$ (with braid index $4$)
gives two $4$-braided fibered patterns and
infinitely many mutually non-conjugate $4$-braided
non-fibered patterns.

Let $\beta_n$ denote the $4$-braid 
$
\sigma_2{\sigma_1}^{-1}{\sigma_2}^{-2n}{\sigma_3}^{-1}
\sigma_2{\sigma_3}^{-1}\sigma_2^{2n}
{\sigma_1}^{-2}$, where $n\ge 0$. Let
$\hat{\beta_n}$ be its closure along an axis $A$.
Each $\hat{\beta_n}$ is (the mirror
image of) the knot $7_6$. 

Let $N(A)$ denote a closed tubular neighborhood of $A$ that is disjoint from $\hat{\beta_n}$. Then $\tilde V={\Bbb S}^3 \setminus int\ N(A)$ is a solid torus, and $(\tilde V, \hat{\beta_n})$ is a pattern of $7_6$. 

\begin{proposition}\label{3componentexample} With $\beta_n$ defined as above:

\begin{enumerate} 
\item \quad The pattern $(\tilde V, \hat{\beta_n})$ is not fibered for any $n> 1$;

\item \quad Patterns $(\tilde V, \hat{\beta_0})$ and $(\tilde V, \hat{\beta_1})$ are fibered.

\end{enumerate}

  \end{proposition}

\begin{remark} Since $\beta_0$ is a homogeneous braid, J. Stallings \cite{stallings2} and L. Rudolph \cite{rudolph} claim that $(\tilde V, \hat{\beta_0})$ is a fibered pattern. However, for completeness we include a proof for this case.

Then encircled closures of $\beta_0$ and $\beta_1$ are links with distinct 1-variable Alexander polynomials, and hence the braids are non-conjugate. 

The braids $\beta_n$, $n\ge 2$, will be seen to be mutually non-conjugate as a consequence of the fact that the genera of the $3$-component links
consisting the closed braids and two copies of axes with
opposite orientations are distinct. Consequently, the patterns
$(\tilde V, \hat{\beta_n})$ are seen to be pairwise distinct.

 \end{remark}

{\it Proof.} For any $n\ge 0$, let $L_n$ denote the $3$-component link consisting
of $\hat{\beta_n}$ and $A\cup A'$, two copies of $A$
with opposite orientations. We will construct a minimal genus Seifert surface
for $L_n$, and show that it is not a fiber surface unless $n=0,1$. We use sutured manifold theory. (For basic facts about sutured manifolds,
see \cite[pp. 8 -- 10 and Appendix A]{Ga1} and \cite[Section 1]{Ga2}.)

\begin{figure}\begin{center}
\includegraphics[height=1.5 in]{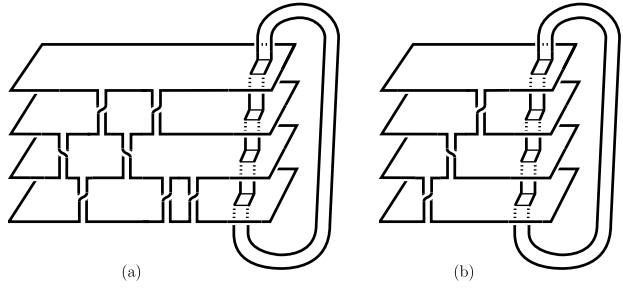}
\caption{Seifert surfaces for $3$-component links}
\label{A1}
\end{center}\end{figure}

Consider the Seifert surface $F_0$ for
$L_0$ as in Figure \ref{A1} (a), obtained by 
smoothing the ribbon intersections of the Seifert surface
for $\hat{\beta_0}$ and the annulus spanned by $A\cup A'$.
Recall that if a Seifert surface has minimal genus or is a fiber, then it remains so after plumbing or deplumbing of Hopf bands. 
Let $F_0'$ be the surface 
as in Figure \ref{A1} (b), obtained
from $F_0$ by deplumbing Hopf bands.

\begin{figure}\begin{center}
\includegraphics[height=2.2 in]{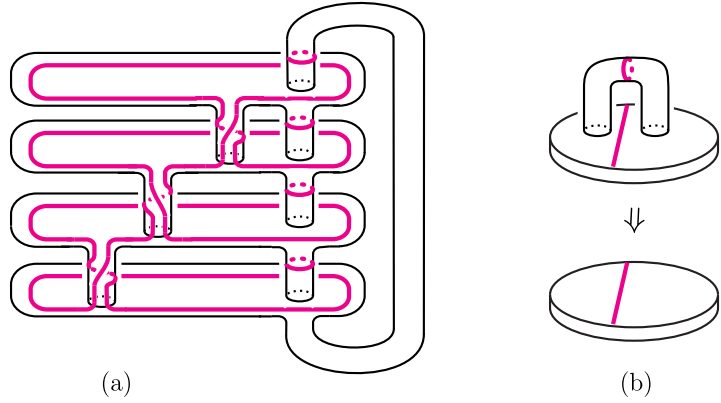}
\caption{A sutured manifold and the $C$-product decomposition}
\label{A2}
\end{center}\end{figure}

Consider the sutured manifold $(M, \gamma)=
(F_0' \times I, \partial F_0' \times I)$,
where $I$ is an interval.
Apply a $C$-procut decomposition at each \lq ribbon hole\rq\
to obtain the sutured manifold as in Figure 
\ref{A2} (a).
The \lq annulus part\rq\ of $F_0$ gives rise to the four
$1$-handles. We can slide each of them 
as in Figure \ref{A2} (b)
and then remove them by $C$-product decompositions.

We now have a sutured manifold which is a $3$-ball with
a single suture. By \cite{Ga1}, $F_0$ is a fiber 
surface (of genus 5) and hence $L_0$ is fibered.
By Theorem \ref{alternate},
 $(\tilde V, \hat{
\beta_0})$ is a fibered pattern.

\begin{figure}\begin{center}
\includegraphics[height=1 in]{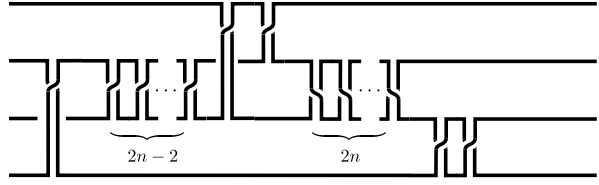}
\caption{Seifert surfaces for closed braids}
\label{A3}
\end{center}\end{figure}

Assume that $n>1$.
Let $R_n$ be a Seifert surface for $\hat{\beta_n}$
(not of minimal genus) as in Figure \ref{A3}.

By superimposing an annulus spanned by $A\cup A'$
and smoothing the ribbon intersections as before,
we obtain a Seifert surface $F_n$ (of genus $2n+3$ if $n>0$) for $L_n$.
Recall that for any fibered link $L$, fiber surfaces are
unique Seifert surfaces of minimal genus.
Therefore, to show that $L_n$ is not a fibered link,
it suffices to show that $F_n$ is not a fiber surface
but that it has minimal genus.

\begin{figure}\begin{center}
\includegraphics[height=1 in]{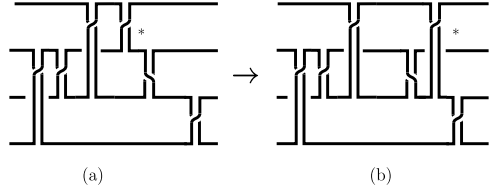}
\caption{Deplumbing a Hopf band}
\label{A4}
\end{center}\end{figure}

Without loss of generality,
we can remove some of parallel half-twisted bands by
deplumbing Hopf bands so that the remaining ones are
as in Figure \ref{A4} (a).
Then we slide a band as in Figure 
\ref{A4} (b) and remove
the band marked $*$ by another Hopf deplumbing.

\begin{figure}\begin{center}
\includegraphics[height=2 in]{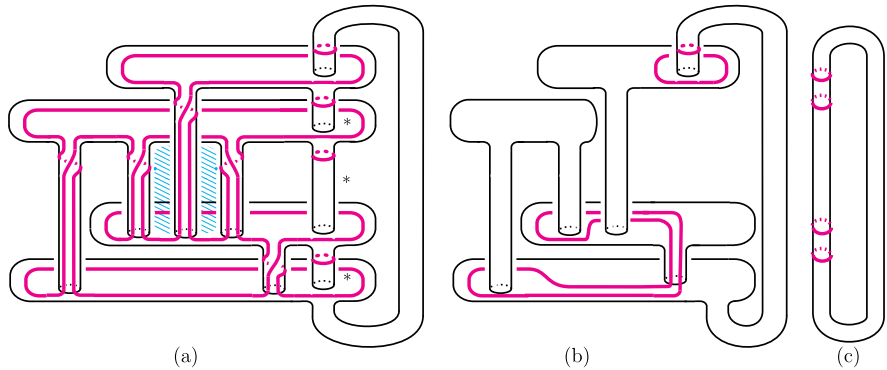}
\caption{$C$-product decompositions}
\label{A5}
\end{center}\end{figure}

Obtain a sutured manifold by thickening the resulting
Seifert surface and apply $C$-product decompositions
at each  cite of a ribbon singularity. The sutured manifold 
appears in Figure \ref{A5} (a).

As before, we can remove three $1$-handles
marked $*$ by $C$-product decompositions.
Then applying another $C$-product decomposition
along the shaded disk (where the dots indicate
the intersections with the sutures),
we obtain a sutured manifold as in Figure 
\ref{A5} (b).

A final $C$-product decomposition yields
a sutured manifold $(M, \gamma)$.
Since it is a solid torus with
four meridional sutures (Figure \ref{A5} (c)), 
its complement
is not a product sutured manifold.

By \cite{Ga1}, $F_n$ is not a fiber surface.
However, we can apply a disk decomposition
to $(M, \gamma)$ so that the result is a $3$-ball 
with a single suture, and hence by \cite{Ga2},
$F_n$ is a minimal
genus Seifert surface. Therefore, $L_n, n>1$ is not
a fibered link.
By Theorem \ref{alternate},
we conclude that $(\tilde V, \hat{\beta_n})$
is not fibered for any $n>1$. 

In the remaining case $n=1$, a sequence of 
$C$- product decompositions as above produces 
a solid torus with two meridional sutures. Consequently,
$F_1$ is  a fiber for $L_1$.  \qed

\section{Fibered patterns of the unknot.}  
We will denote by ${\cal P}_n$ the collection of all fibered patterns $(\tilde V, \tilde k)$ such that $\tilde k$ is a trivial knot with winding number $n$ in $\tilde V$. The following theorem classifies the members of ${\cal P}_n$.

\begin{theorem}\label{classify} Let $(\tilde V, \tilde k)$ be a pattern such that $\tilde k$ is the trivial knot with winding number in $\tilde V$ equal to $n$. Then $(\tilde V, \tilde k)$ is fibered if and only if the meridian $\tilde m$ of $\tilde V$ is the closure of an $n$-braid $B$ with axis $\tilde k$. The mapping $(\tilde V, \tilde k) \mapsto [B]$ is a bijection between ${\cal P}_n$  and the set of conjugacy classes of the $n$-braid group $B_n$. \end{theorem}

\noindent For the definition and elementary properties of $n$-braids, the reader is referred to \cite{birman}.

{\it Proof.} We use notation established in section 2. Assume that $(\tilde V, \tilde k)$ is a member of ${\cal P}_n$. Then there exists a fibration $Y(\tilde k) \to \S^1$ inducing the homomorphism $\phi\circ \pi$ on fundamental groups. Let $S$ be a fiber $p^{-1}(\theta),\ \theta \in \S^1$. Since $p$ restricts to a fibration of 
$\partial Y(\tilde k)$, it follows that $\partial S$ consists of a longitude of $\tilde k$ together with longitudes $\tilde \ell_1, \ldots, \tilde \ell_n$ of $\tilde V$. We can extend $p$ to a fibration $\bar p: X(\tilde k) \to \S^1$ such that $\bar p^{-1}(\theta)$ is the ``capped off" fiber $S \cup_{\tilde \ell_1} \D^2 \cup_{\bar \ell_2} \cdots \cup_{\bar \ell_n} \D^2$, since $\tilde V$ is unknotted. Assume now that $\tilde k$ is a trivial knot. Then $\bar p^{-1}(\theta)$ must be a disk. Moreover, $X(\tilde k)$ is isotopic in $\S^3$ to a standardly embedded solid torus $\S^1 \times \D^2$.  Since $\tilde m$ intersects each $\theta \times \D^2\ = (\bar p^{-1}(\theta))$ in exactly $n$ points, $\tilde m$ is the closure of an $n$-braid $B$ with axis $\tilde k$. The particular embedding of $\tilde m$  in $\S^1 \times \D^2$ that we obtain in this way is independent up to isotopy of the choices that we make. Hence by Theorem 1 of \cite{morton}, the conjugacy class of $B$ in the $n$-braid group is well defined. 

In order to prove that the mapping $(\tilde V, \tilde k) \mapsto [B]$ is a bijection, it suffices to produce the inverse mapping. The method that we use is essentially the reverse of the above, and it has been used elsewhere (cf. \cite{morton}). Assume that $B$ is an $n$-braid that closes to a trivial knot $\tilde m$. Denote the axis of $\tilde m$ by $\tilde k$. Let $N(\tilde m)$ and $N(\tilde k)$ be disjoint tubular neighborhoods  in $\S^3$ of $\tilde m$ and $\tilde k$, respectively, and consider the space $Y= \S^3 \setminus int\ (N(\tilde m) \cup N(\tilde k))$. Regarded as the exterior of the closed $n$-braid $\tilde m$ in the solid torus $\S^3 \setminus int\ N(\tilde k)$, the space $Y$ fibers over $\S^1$; the fiber is a disk with $n$ holes and boundary consisting of a meridian  of $\S^3 \setminus int\ N(\tilde k)$ together with $n$ meridians of $N(\tilde m)$. On the other hand, $Y$ can be viewed as the exterior of the trivial knot $\tilde k$ in the unknotted solid torus $\tilde V = \S^3 \setminus int\ N(\tilde m)$. From this vantage, the boundary of the fiber consists of a longitude of $\tilde k$ together with $n$ longitudes of $\tilde V$, each also a meridian of $\S^3 \setminus int \tilde V$. It follows that the pattern $(\tilde V, \tilde k)$ is fibered. Notice that $\tilde k$ has winding number in $\tilde V$ equal to $n$. Also, $(\tilde V, \tilde k)$ depends only on the conjugacy class of $B$. Clearly, $B \mapsto (\tilde V, \tilde k)$ is the desired inverse mapping. 
\qed

If $\tilde m$ is the the closure of a $1$-braid with unknotted axis $\tilde k$, then $\tilde k \cup \tilde m$ is a Hopf link. We obtain the following corollary. In effect, it asserts that there is no fibered satellite knot with unknotted pattern knot $\tilde k$ having winding number $1$ in the solid torus $\tilde V$. 

\begin{corollary} \label{unknotted} Let $\tilde k$ be a trivial knot embedded in an unknotted solid torus $\tilde V$ with winding number equal to $1$. Then $(\tilde V, \tilde k)$ is fibered if and only if $(\tilde V, \tilde k) = (\S^1 \times \D^2, \S^1 \times 0)$. 
\end{corollary}

\begin{remark} There exist nontrivial fibered patterns with winding number $1$. Such an example $(\tilde{V}, \tilde{k})$ can be constructed from the closure of the $3$-braid $\sigma_2^2 \sigma_1^2 \sigma_2$ in $\tilde{V}$ by adding a short band with a full twist between the
1st and 2nd strings, away from the crossings.
Here $\tilde{k}$ naturally spans a fiber surface $F$ consisting of a checkerboard-colorable surface together with a plumbed Hopf band. A fiber $S$ for the $3$-component link of Theorem \ref{alternate} is obtained by adding an annulus to $F$. By Theorem \ref{alternate} and arguments similar to earlier ones, we see that $(\tilde{V}, \tilde{k})$ is a fibered pattern. 

\end{remark}

\begin{figure}
\begin{center}
\includegraphics[height=1.6in]{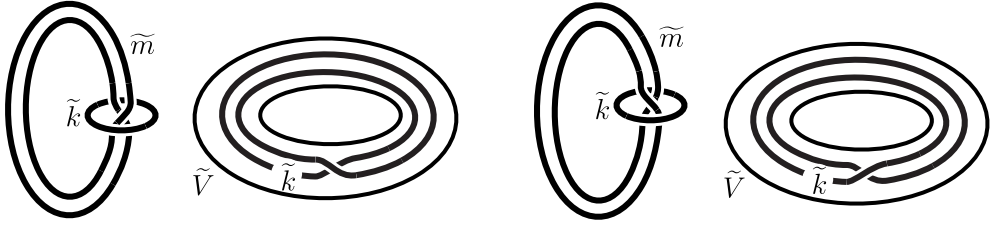}
\caption{Fibered patterns of unknot with winding number 2.}
\label{2braids}
\end{center}
\end{figure}

By reasoning similar to that of Corollary \ref{unknotted} we find that there exist only two fibered patterns $(\tilde V, \tilde k)$ such that $\tilde k$ is a trivial knot with winding number in $\tilde V$ equal to $2$. These correspond to the conjugacy classes of the $2$-braids $\sigma_1$ and $\sigma_1^{-1}$. 
See \Fig. \ref{2braids}. 

For $n=3$, there are only three fibered patterns, corresponding to the conjugacy classes of $3$-braids $\sigma_1 \sigma_2,\ \sigma_1^{-1} \sigma_2$ and $\sigma_1^{-1}\sigma_2^{-1}$ \cite{mp}. (The braid $\sigma_1 \sigma_2^{-1}$ does not appear in the list since it is conjugate to 
$\sigma_1^{-1} \sigma_2$ in $B_3$.)

If $n$ is greater than 3, the situation is more complicated. We note that  
$\sigma_1 \sigma_2^{2i+1} \sigma_3 \sigma_2^{-2i},\ i\ge 0$ represent infinitely many pairwise non-conjugate $4$-braids
\cite{morton}. Axes $\tilde m$ of non-conjugate braids form distinct fibered patterns of the unknot by Theorem \ref{classify}.  Hence ${\cal P}_4$ is infinite.

\begin{corollary} 
Any knot $\tilde k$ can be embedded in an unknotted solid torus $\tilde V$ with arbitrary winding number $n$ in such a way that no satellite knot with pattern $(\tilde V, \tilde k)$ is fibered. 
\end{corollary}

{\it Proof.} By Theorem  \ref{criterion1} it suffices to prove that $\tilde k$ can be embedded in $\tilde V$ with arbitrary winding number $n$ such that $(\tilde V, \tilde k)$ is not fibered. We may assume that $n$ is positive (see remarks preceding Theorem \ref{criterion1}). We will prove the corollary assuming further that $\tilde k$ is the trivial knot. The general result will then follow by trying an arbitrary ``local knot" in $\tilde k$ (that is, forming the knot product of $\tilde k$ and an arbitrary knot, and doing so within a $3$-ball in $\tilde V$), since a product of knots is fibered if and only if each of the factors is.

\begin{figure}
\centering
 \includegraphics [height=2in]
 {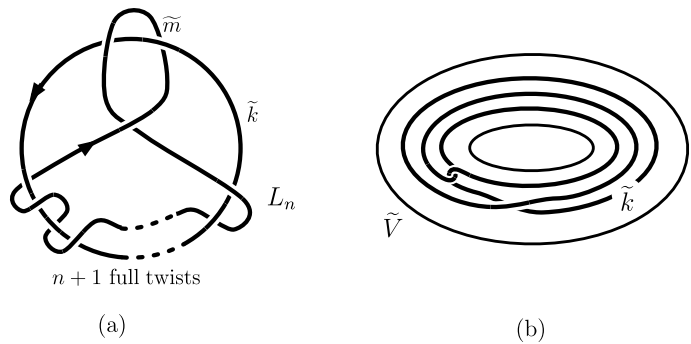}
\caption{(a) The link $L_n$, (b) Pattern when $n=1$.}
\label{link}
\end{figure}

For $n>0$, consider the oriented link $L_n$ depicted in \Fig. \ref{link} (a),
consisting of two unknotted components $\tilde k$ and $\tilde m$.

Let $N(\tilde m)$ be a tubular neighborhood in $\S^3$ of $\tilde m$ that is disjoint from $\tilde k$. Then $\S^3 \setminus int\ N(\tilde m)$ is an unknotted solid torus $\tilde V$ containing the trivial knot $\tilde k$ with winding number $n$. 
(When $n=1$, the pattern $(\tilde V, \tilde k)$ is the configuration in \Fig. \ref{link} (b) 
(cf. \cite{mazur})).  Assume that $(\tilde V, \tilde k)$ is fibered. By Theorem \ref{classify}, the link $L_n$ can be isotoped in such a way that $\tilde m$ appears as the closure of  an $n$-braid with axis $\tilde k$. Consequently, $L_n$ is the closure of an $(n+1)$-braid. A theorem of J. Franks and R.F. Williams \cite{fw} (also proved by H.R. Morton) implies that the braid index of 
$L_n$ is at least
$\frac{1}{2}(e_{max} - e_{min})+1$, where $e_{max}$ and $e_{min}$ are respectively the largest and smallest exponents of $x$ in the $2$-variable Jones polynomial $j_{L_n}(x,y)$ of $L_n$. 
An easy calculation reveals that $j_{L_n}(x, y) = -x^{-2}j_{L_{n-1}}(x,y) - x^{-1} j_{3_1}(x,y)$, where $3_1$ denotes the left-hand trefoil knot, and from this it is not difficult to verify that ${1\over 2}(e_{max} - e_{min})+1 = n+3$, a contradiction. Hence $(\tilde V, \tilde k)$ is not fibered. 
\qed

\begin{remark} The orientation of the components of $L_n$ is critical. If the orientation of just one component is reversed, then the braid index of $L_n$ becomes 3, for any nonzero $n$. 

\end{remark}

\end{document}